\documentclass[10pt]{article}
\usepackage{amsmath, amssymb, amscd, amsthm}

\newcommand{\xra}[1]{\ensuremath{\xrightarrow{#1}}}
\newcommand{\B}[1]{\ensuremath{\mathbb{#1}}}
\newcommand{\C}[1]{\ensuremath{\mathcal{#1}}}
\newcommand{\G}[1]{\ensuremath{\mathfrak{#1}}}
\newcommand{\eotimes}[1]{\ensuremath{\underset{#1}{\otimes}}}
\newcommand{\fotimes}[1]{\ensuremath{\underset{#1}{\Box}}}

\newtheorem{thm}{Theorem}[section]
\newtheorem{cor}[thm]{Corollary}
\newtheorem{lem}[thm]{Lemma}

\theoremstyle{definition}
\newtheorem{defn}[thm]{Definition}
\newtheorem{rem}[thm]{Remark}
\newtheorem{exm}[thm]{Example}

\numberwithin{equation}{section}

\title{Bialgebra Cyclic Homology with Coefficients\\ Part II}
\author{Atabey Kaygun}
\date{}
\begin{document}
\maketitle

\section{Introduction}

This is the second part of the article \cite{Kaygun:BialgebraCyclic}.
In the first paper we developed a cyclic homology theory for $B$--module
coalgebras with coefficients in stable $B$--module/comodules where $B$
was just a bialgebra.  The construction we gave for the cyclic homology
theory for $B$--module coalgebras used mainly the coalgebra structure on
$B$.  In the first part of this paper, we present the dual picture.
Namely, a cyclic homology theory for $B$--comodule algebras with
coefficients in a stable $B$--module/comodule where $B$ is just a
bialgebra.  Our theory is an extension of the theory developed in
\cite{Khalkhali:HopfCyclicHomology} by lifting two restrictions: (i) our
theory uses bialgebras as opposed to Hopf algebras (ii) the coefficient
module/comodules are just stable as opposed to stable
anti-Yetter-Drinfeld.  In the second part of this paper, we recover the
main result of \cite{Khalkhali:CyclicDuality}.  Namely, these two cyclic
theories are dual in the sense of (co)cyclic objects, whenever the input
pair $(H,X)$ has the property that $H$ is a Hopf algebra and $X$ is a
stable anti-Yetter-Drinfeld module.

The plan of this paper is as follows.  In Section~\ref{Notation}, we set
up the notation and overall assumptions we make.  In
Section~\ref{OldTheory}, we develop a cyclic theory for a pair $(H,X)$
where $H$ is a Hopf algebra and $X$ is a stable anti-Yetter-Drinfeld
module by using the algebra structure of $H$ and $H$--module/comodule
structure of $X$.  In Section~\ref{NewTheory}, we show how one can
extend this theory to bialgebra comodule algebras and stable bialgebra
modules.  In Section~\ref{Duality}, we show that the cyclic theory we
developed in \cite{Kaygun:BialgebraCyclic} and cyclic theory defined in
this paper are dual in the sense of (co)cyclic objects whenever the
underlying bialgebra is a Hopf algebra and the stable coefficient
module/comodule is also anti-Yetter-Drinfeld.  In
Section~\ref{Computations} we perform several calculations to illustrate
the effectiveness of our definition of bialgebra cyclic homology for
several Hopf algebras: the group ring of a discrete group $G$, the
enveloping algebra of a Lie algebra $\G{g}$, quantum deformation of an
arbitrary semi-simple Lie algebra $\G{g}$, and finally $\C{H}(N)$, the
Hopf algebra of foliations of codimension $N$.

\section{Notation and conventions}\label{Notation}

\noindent We assume $k$ is a field of an arbitrary characteristic and
$H$ is a Hopf algebra over $k$.

\noindent Whenever we refer an object ``simplicial'' or
``cosimplicial,'' the reader should read as ``pre-simplicial'' and
``pre-cosimplicial'' meaning that do not consider (co)degeneracy
morphisms as a part of the (co)simplicial data.

\noindent A simplicial $X_*$ module is called para-(co)cyclic iff it is
almost a (co)cyclic module, in that it satisfies all conditions for a
(co)cyclic module except that the action of $\tau_n$ on each $X_n$ need
not to be of order $n+1$, for any $n\geq 0$.

\noindent A (para-)cyclic module $\C{Z}_*$ is called a (para-)cyclic
$H$--comodule iff all structure morphisms are $H$--comodule morphisms.

\noindent The tensor product over an algebra $A$ is denoted by
$\eotimes{A}$, and a cotensor product over a coalgebra $C$ is denoted
by $\fotimes{C}$.  Recall that if $X\xra{\rho_X}X\otimes C$ and
$Y\xra{\rho_Y}C\otimes Y$ are two $C$--comodules (right and left
respectively), then $X\fotimes{C}Y$ is defined as $ker((\rho_X\otimes
id_Y)-(id_X\otimes \rho_Y))$. 

\noindent For a coalgebra $(C,\Delta)$, we use Sweedler's notation and
denote $\Delta(c)$ by $\sum_c c_{(1)}\otimes c_{(2)}$, and most of the
time, we even drop the summation sign.  Similarly, for a left
$C$--comodule $X\xra{\rho_X}C\otimes X$, we use $\rho_X(x) =
x_{(-1)}\otimes x_{(0)}$ for the coaction morphism.  On the other hand,
for a right $C$--comodule $Y\xra{\rho_Y}Y\otimes C$ we use
$\rho_Y(y)=(y_{(0)}\otimes y_{(1)})$.

\noindent Given a counital bialgebra $(B,\cdot,\Delta,\epsilon)$ and a
right $B$--comodule $M\xra{\rho_M}M\otimes B$, the comodule of
$B$--invariants of $M$ which is
\[ \left\{m|\ \rho_M(m)= (m_{(0)}\otimes m_{(1)}) = (m\otimes\B{I})\right\} \] 
is denoted by $M^B$.

\noindent We also need to use complexes $\B{T}_*(B,X)$ and
$\B{CM}_*(B,X)$ we defined in \cite{Kaygun:BialgebraCyclic}.  In order
to distinguish these complexes from the similar complexes we define in
this paper, we use the notation $\B{T}^c_*(B,X)$ and $\B{CM}^c_*(B,X)$.

\section{The Connes--Moscovici cyclic homology}\label{OldTheory}

\begin{defn}
Let $A$ be an algebra over $k$, let $X$ be a left $A$--module and let
$Y$ be a right $A$--module.  The bar complex associated to the algebra
$A$ with coefficients in the $A$--modules $A$ and $Y$ is the simplicial
$k$--module $B^a_*(Y,A,X)=\{Y\otimes A^{\otimes n}\otimes X\}_{n\geq 0}$
with the following face morphisms:
\begin{equation*}
d_j(y\otimes a^1\otimes\cdots\otimes a^n\otimes x)
 = \begin{cases}
   (ya^1\otimes\cdots\otimes x)   & \text{ if } j=0\\
   (y\otimes\cdots\otimes a^j a^{j+1}\otimes\cdots\otimes x)
                                  & \text{ if } 0<j<n\\
   (y\otimes\cdots\otimes  h^nx)  & \text{ if } j=n          
   \end{cases}
\end{equation*}
for any $(y\otimes a^1\otimes\cdots\otimes a^n\otimes x)$ from
$B_n(Y,A,X)$.
\end{defn}

\begin{defn}
Let $\B{T}^a_*(H,X)=\{H^{\otimes n+1}\otimes X\}_{n\geq 0}$ where $X$ is a
$H$--module.  Define
\begin{align}
\partial_j(h^0\otimes\cdots\otimes h^n\otimes x)
 = & \begin{cases}
     (\cdots\otimes h^jh^{j+1}\otimes\cdots\otimes x)
        & \text{ if } 0\leq j< n\\
     \left(h^n_{(1)}h^0\otimes\cdots\otimes h^{n-1}\otimes 
           h^n_{(2)}x\right)
        & \text{ if } j=n
     \end{cases}
\end{align}
for $0\leq j\leq n$ and for any $(h^0\otimes\cdots\otimes h^n\otimes x)$
from $\B{T}^a_n(H,X)$.  Then $\B{T}^a_*(H,X)$ is a simplicial
$k$--module.
\end{defn}

\begin{lem}\label{adjoint}
Let $M$ be a $H$--bimodule.  Then there is a right $H$--module structure
on $M$ defined as $m\cdot ad_h := S^{-1}(h_{(1)})mh_{(2)}$ for any $h\in H$
and $m\in M$.  This action is called the right adjoint action and the
module is denoted by $ad(M)$.
\end{lem}

\begin{thm}\label{bar_complex}
Let $X$ be an arbitrary $H$--module.  There is an isomorphism of
simplicial $k$--modules of the form
$\B{T}^a_*(H,X)\xra{\Phi_*}B^a_*(ad(H),H,X)$ where
\begin{align*}
\Phi_n(h^0\otimes h^1\otimes\cdots\otimes h^n\otimes x)
 = (h^1_{(1)}\cdots h^n_{(1)}h^0\otimes h^1_{(2)}\otimes\cdots\otimes
    h^n_{(2)}\otimes x)
\end{align*}
for any $(h^0\otimes h^1\otimes\cdots\otimes h^n\otimes x)$ from
$\B{T}^a_n(H,X)$.
\end{thm}

\begin{proof}
The inverse is given by
\begin{equation*}
\Phi^{-1}_n(h^0\otimes h^1\otimes\cdots\otimes h^n\otimes x)
 = \left(S^{-1}(h^1_{(1)}\cdots h^n_{(1)})h^0\otimes
         h^1_{(2)}\otimes\cdots\otimes h^n_{(2)}\otimes x\right)
\end{equation*}
since
\begin{align*}
\Phi_n\Phi^{-1}_n(h^0\otimes h^1\otimes\cdots\otimes h^n\otimes x)
 = & \Phi_n\left(S^{-1}(h^1_{(1)}\cdots h^n_{(1)})h^0\otimes
                 h^1_{(2)}\otimes\cdots\otimes h^n_{(2)}
                 \otimes x\right)\\
 = & \left(h^1_{(2)(1)}\cdots h^n_{(2)(1)}
           S^{-1}(h^1_{(1)}\cdots h^n_{(1)})h^0\otimes
           h^1_{(2)(2)}\otimes\cdots\otimes h^n_{(2)(2)}
           \otimes x\right)\\
 = & (h^0\otimes h^1\otimes\cdots\otimes h^n\otimes x)
\end{align*}
and similarly
\begin{align*}
\Phi^{-1}_n\Phi_n(h^0\otimes h^1\otimes\cdots\otimes h^n\otimes x)
 = & \Phi^{-1}_n\left(h^1_{(1)}\cdots h^n_{(1)}h^0\otimes 
                      h^1_{(2)}\cdots\otimes h^n_{(2)}
                      \otimes x\right)\\
 = & \left(S^{-1}(h^1_{(2)(1)}\cdots h^n_{(2)(1)})
           h^1_{(1)}\cdots h^n_{(1)}h^0\otimes 
           h^1_{(2)(2)}\cdots\otimes h^n_{(2)(2)}\otimes x\right)\\
 = & (h^0\otimes h^1\otimes\cdots\otimes h^n\otimes x)
\end{align*}
as we wanted to show.  Now consider
\begin{align*}
\Phi_{n-1}\partial_0(h^0\otimes h^1\otimes\cdots\otimes h^n\otimes x)
 = & \Phi_{n-1}(h^0 h^1\otimes\cdots\otimes h^n\otimes x)\\
 = & \left(h^2_{(1)}\cdots h^n_{(1)}h^0h^1\otimes 
           h^2_{(2)}\cdots\otimes h^n_{(2)}\otimes x\right)\\
 = & \left((h^1_{(1)}\cdots h^n_{(1)}h^0)\cdot ad_{h^1_{(2)}}\otimes 
            h^2_{(2)}\cdots\otimes h^n_{(2)}\otimes x\right)\\
 = & d_0\Phi_n(h^0\otimes h^1\otimes\cdots\otimes h^n\otimes x)
\end{align*}
and for $0<j<n$
\begin{align*}
\Phi_{n-1}\partial_j(h^0\otimes h^1\otimes\cdots\otimes h^n\otimes x)
 = & \Phi_{n-1}(h^0\otimes\cdots\otimes h^jh^{j+1}\otimes\cdots
                \otimes x)\\
 = & \left(h^1_{(1)}\cdots h^n_{(1)}h^0\otimes\cdots\otimes
      h^j_{(2)}h^{j+1}_{(2)}\otimes\cdots\otimes x\right)\\
 = & d_j\Phi_n(h^0\otimes\cdots\otimes h^n\otimes x)
\end{align*}
Finally for $j=n$,
\begin{align*}
\Phi^{-1}_{n-1}d_n\Phi_n(h^0\otimes h^1\otimes\cdots
        \otimes h^n\otimes x)
 = & \Phi^{-1}_{n-1}\left(h^1_{(1)}\cdots h^n_{(1)}h^0\otimes h^1_{(2)}
        \otimes\cdots\otimes h^n_{(2)}x\right)\\
 = & \left(h^n_{(1)}h^0\otimes h^1\otimes\cdots\otimes h^{n-1}\otimes
          h^n_{(2)}x\right)
\end{align*}
as we wanted to prove.
\end{proof}

\begin{defn}
Let $X$ be an arbitrary $H$--module.  Define a para-cyclic structure on
$\B{T}^a_*(H,X)$ by letting
\begin{align}
\partial_0(h^0\otimes\cdots\otimes h^n\otimes x)
 = &  (h^0 h^1\otimes\cdots\otimes h^n\otimes x)\\
\tau_n(h^0\otimes\cdots\otimes h^n\otimes x)
 = & \left(h^n_{(1)}\otimes h^0\otimes\cdots\otimes h^{n-1}\otimes 
       h^n_{(2)}x\right)\label{cyclic_definition}\\
\tau^{-1}_n(h^0\otimes\cdots\otimes h^n\otimes x)
 = & \left(h^1\otimes\cdots\otimes h^n\otimes 
        h^0_{(1)}\otimes S(h^0_{(2)})x\right)\\
\partial_j (h^0\otimes\cdots\otimes h^n\otimes x)
 = & \tau_{n-1}^j\partial_0\tau_n^{-j}
       (h^0\otimes\cdots\otimes h^n\otimes x)
\end{align}
for $0\leq j\leq n$.
\end{defn}

\begin{lem}~\label{face}
Let $n\geq 0$ be arbitrary and let $(h^0\otimes\cdots\otimes h^n\otimes
x)$ be from $\B{T}^a_n(H,X)$.  Then
\begin{align*}
\tau_{n-1}^{-n}\partial_0\tau_n^{n+1}
     (h^0\otimes\cdots\otimes h^n\otimes x)
 = & \partial_0(h^0\otimes\cdots\otimes h^n\otimes x)
\end{align*}
\end{lem}

\begin{proof}
Consider 
\begin{align*}
\tau_n^{n+1}(h^0\otimes\cdots\otimes h^n\otimes x)
 = & \left(h^0_{(1)}\otimes\cdots\otimes h^n_{(1)}\otimes
           h^0_{(2)}\cdots h^n_{(2)}x\right)
\end{align*}
then
\begin{align*}
\partial_0\tau_n^{n+1}(h^0\otimes\cdots\otimes h^n\otimes x)
 = & \left(h^0_{(1)}h^1_{(1)}\otimes\cdots\otimes h^n_{(1)}\otimes
           h^0_{(2)}\cdots h^n_{(2)}x\right)
\end{align*}
and finally
\begin{align*}
\tau_n^{-n}\partial_0\tau_n^{n+1}(h^0\otimes\cdots\otimes h^n\otimes x)
 = & \left(h^0_{(1)(1)}h^1_{(1)(1)}\otimes\cdots\otimes h^n_{(1)(1)}\otimes
           S(h^n_{(1)(2)})\cdots S(h^0_{(1)(2)})
             h^0_{(2)}\cdots h^n_{(2)}x\right)\\
 = & (h^0h^1\otimes\cdots\otimes h^n\otimes x)\\
 = & \partial_0(h^0\otimes\cdots\otimes h^n\otimes x)
\end{align*}
for any $n\geq 0$ and for any $(h^0\otimes\cdots\otimes h^n\otimes x)$
from $\B{T}^a_n(H,X)$.
\end{proof}

\begin{cor}
The face morphisms in $\B{T}^a_*(H,X)$ are defined as
\begin{align}
\partial_j({\bf h}\otimes x)
 = \tau_{n-1}^{-n+j}\partial_0\tau_n^{n+1-j}({\bf h}\otimes x)
 = \begin{cases}
     (\cdots\otimes h^jh^{j+1}\otimes\cdots\otimes x)
        & \text{ if } 0\leq j< n\\
     \left(h^n_{(1)}h^0\otimes\cdots\otimes h^{n-1}\otimes 
           h^n_{(2)}x\right)
        & \text{ if } j=n
     \end{cases}\label{simplicial_definition}
\end{align}
for any $n\geq 0$ and $({\bf h}\otimes x)$ from $\B{T}^a_n(H,X)$.
\end{cor}

\begin{defn}
Let $X$ be an arbitrary $H$--module/comodule.  Define a graded
$k$--module by $\B{CM}^a_*(H,X)=\{H^{\otimes n}\otimes X\}_{n\geq 0}$ and
a pair of graded $k$--module morphisms
$\B{CM}^a_*(H,X)\xra{p_*}\B{T}^a_*(H,X)$ and
$\B{T}^a_*(H,X)\xra{i_*}\B{CM}^a_*(H,X)$ by
\begin{align*}
i_n(h^0\otimes\cdots\otimes h^n\otimes x) 
 = & \left(h^0\otimes\cdots\otimes h^{n-1}\otimes h^n x\right)\\
p_n(h^1\otimes\cdots\otimes h^n\otimes x)
 = & \begin{cases}
     (x_{(1)}\otimes x_{(0)})
           & \text{ if } n=0\\
     \left(h^1_{(1)}\otimes\cdots\otimes h^n_{(1)}\otimes
           x_{(-1)}S^{-1}(h^1_{(3)}\cdots h^n_{(3)})\otimes
           h^1_{(2)}\cdots h^n_{(2)}x_{(0)}\right)
           & \text{ if } n>0
     \end{cases}
\end{align*}
\end{defn}

\begin{defn}
Let $H$ be a Hopf algebra.  Then a $H$--module/comodule $X$ is called
$m$-stable if \[ S^m(x_{(-1)})x_{(0)}=y \] for all $x\in X$.  If $X$ is
both $1$--stable and $0$--stable, we call it stable.
\end{defn}

\begin{defn}
Let $H$ be a Hopf algebra.  Then a $H$--module/comodule is called
anti-Yetter-Drinfeld (aYD) module iff
\begin{align*}
(hx)_{(-1)}\otimes (hx)_{(0)} = h_{(1)}x_{(-1)}S^{-1}(h_{(3)})\otimes
  h_{(2)}x_{(0)} 
\end{align*}
for any $x\in X$ and $h\in H$.  
\end{defn}

\begin{lem}
Assume $X$ is an anti-Yetter-Drinfeld module.  Then $X$ is $0$-stable
iff $X$ is $1$-stable.
\end{lem}

\begin{proof}
Assume $x_{(-1)}x_{(0)}=x$ for any $x\in X$.  Let $y=S(x_{(-1)})x_{(0)}$
and consider
\begin{align*}
y = & y_{(-1)}y_{(0)} \\
  = & S(x_{(-1)(3)})x_{(0)(-1)}x_{(-1)(1)}S(x_{(-1)(2)})x_{(0)(0)} \\
  = & S(x_{(-2)})x_{(-1)}x_{(-4)}S(x_{(-3)})x_{(0)}\\
  = & x
\end{align*}
The proof for the other direction is similar.
\end{proof}

\begin{rem}\label{monomorphism}
Notice that if we assume $X$ is $0$--stable $H$--module/comodule,
i.e. $x=x_{(-1)}x_{(0)}$ for any $x\in X$, then $i_*p_*=id_*$.  This
implies $p_*$ is a monomorphism of graded $k$--modules.
\end{rem}

\begin{thm}
Assume $X$ is a stable anti-Yetter-Drinfeld module.  Then there is a
para-cyclic structure on $\B{CM}^a_*(H,X)$ such that $\B{CM}^a_*(H,X)
\xra{p_*}\B{T}^a_*(H,X)$ is a morphism of para-cyclic modules.
\end{thm}
\begin{proof}
Define a morphism $d_0$ of degree $-1$ on $\B{CM}^a_*(H,X)$ by letting
\begin{align*}
d_0(h^1\otimes\cdots\otimes h^n\otimes x)
 = \begin{cases}
   h^1x  & \text{ if } n=1\\
   (h^1h^2\otimes\cdots\otimes x)
         & \text{ if } n>1
   \end{cases}
\end{align*}
and observe that if we assume that $X$ is an anti-Yetter-Drinfeld
module, we get
\begin{align*}
p_{n-1}d_0 & (h^1\otimes\cdots\otimes h^n\otimes x)\\
 = & \begin{cases}
     p_0(h^1x)   & \text{ if } n=1\\
     p_{n-1}(h^1h^2\otimes\cdots\otimes x)
                 & \text{ if } n>1
     \end{cases}\\
 = & \begin{cases}
     \left(h^1_{(1)}x_{(-1)}S^{-1}(h^1_{(3)})\otimes
           h^1_{(2)}x_{(0)}\right)   
                 & \text{ if } n=1\\ 
     \left(h^1_{(1)}h^2_{(1)}\otimes\cdots\otimes h^n_{(1)}\otimes
           x_{(-1)}S^{-1}(h^1_{(3)}h^2_{(3)}\cdots h^n_{(3)})\otimes
           h^1_{(2)}h^2_{(2)}\cdots h^n_{(2)}x_{(0)}\right)
                 & \text{ if } n>1
     \end{cases}\\
 = & \partial_0p_n(h^1\otimes\cdots\otimes h^n\otimes x)
\end{align*}
Now let $t_* = i_*\tau_*p_*$.  We need to show that $p_*t_* =
p_*i_*\tau_* p_* = \tau_*p_* $.  So, consider
\begin{align*}
\tau_np_n & (h^1\otimes\cdots\otimes h^n)\nonumber\\
 = & \tau_n\left(h^1_{(1)}\otimes\cdots\otimes h^n_{(1)}\otimes
           x_{(-1)}S^{-1}(h^1_{(3)}\cdots h^n_{(3)})\otimes
           h^1_{(2)}\cdots h^n_{(2)}x_{(0)}\right)\\
 = & \left(x_{(-1)(1)}S^{-1}(h^1_{(3)(2)}\cdots h^n_{(3)(2)})
           \otimes h^1_{(1)}\otimes\cdots\otimes h^n_{(1)}\otimes
           x_{(-1)(2)}S^{-1}(h^1_{(3)(1)}\cdots h^n_{(3)(1)})
           h^1_{(2)}\cdots h^n_{(2)}x_{(0)}\right)\\
 = & \left(x_{(-1)}S^{-1}(h^1_{(2)}\cdots h^n_{(2)})\otimes 
           h^1_{(1)}\otimes\cdots\otimes h^n_{(1)}\otimes
           x_{(0)}\right)
\end{align*}
Assume $n\geq 1$ and consider also
\begin{align*}
t^{-1}_n & (h^1\otimes\cdots\otimes h^n\otimes x)\nonumber\\
 = & i_n\tau_n^{-1}p_n(h^1\otimes\cdots\otimes h^n\otimes x)\\
 = & i_n\tau_n^{-1}\left(h^1_{(1)}\otimes\cdots\otimes h^n_{(1)}
        \otimes x_{(-1)}S^{-1}(h^1_{(3)}\cdots h^n_{(3)})
        \otimes h^1_{(2)}\cdots h^n_{(2)}x_{(0)}\right)\\
 = & i_n\left(h^2_{(1)}\otimes\cdots\otimes h^n_{(1)}
        \otimes x_{(-1)}S^{-1}(h^1_{(3)}\cdots h^n_{(3)})
        \otimes h^1_{(1)(1)}\otimes S(h^1_{(1)(2)})
        h^1_{(2)}\cdots h^n_{(2)}x_{(0)}\right)\\
 = & \epsilon(h^1_{(1)})\left(h^2_{(1)}\otimes\cdots\otimes h^n_{(1)}
        \otimes x_{(-1)}S^{-1}(h^1_{(3)}\cdots h^n_{(3)})
        \otimes h^1_{(2)}h^2_{(2)}\cdots h^n_{(2)}x_{(0)}\right)
\end{align*}
Now, using these identities we consider
\begin{align*}
\tau_np_nt_n^{-1} & (h^1\otimes\cdots\otimes h^n\otimes x)\nonumber\\
 = & \epsilon(h^1_{(1)})\tau_np_n\left(h^2_{(1)}\otimes\cdots\otimes h^n_{(1)}
        \otimes x_{(-1)}S^{-1}(z_{(3)})
        \otimes z_{(2)}x_{(0)}\right)\\
 = & \epsilon(h^1_{(1)})\left(h^1_{(2)(1)}\cdots h^n_{(2)(1)}
        x_{(0)(-1)}\right.\\
   & \hspace{1.5cm}S^{-1}(h^1_{(2)(3)}\cdots h^n_{(2)(3)})
        S^{-1}\left(h^2_{(1)(2)}\cdots h^n_{(1)(2)}x_{(-1)(2)}
        S^{-1}(h^1_{(3)(1)}\cdots h^n_{(3)(1)})\right)\\
   & \left.\hspace{1.5cm}\otimes
        h^2_{(1)(1)}\otimes\cdots\otimes h^n_{(1)(1)}\otimes
        x_{(-1)(1)}S^{-1}(h^1_{(3)(2)}\cdots h^n_{(3)(2)})
        \otimes h^1_{(2)(2)}\cdots h^n_{(2)(2)}x_{(0)(0)}
     \right)\\
 = & \epsilon(h^1_{(1)})\left(h^1_{(2)}\otimes
        h^2_{(1)}\otimes\cdots\otimes h^n_{(1)}\otimes
        x_{(-1)}S^{-1}(h^1_{(3)}\cdots h^n_{(3)})\otimes 
        h^1_{(2)}\cdots h^n_{(2)}x_{(0)}
     \right)\\
 = & p_n(h^1\otimes\cdots\otimes h^n\otimes x)
\end{align*}
as we wanted to show.  This finishes the proof that $\tau_*p_*=p_*t_*$.
Define a cyclic structure on $\B{CM}^a_*(H,Y)$ by letting
\begin{align*}
d_j = t_{n-1}^jd_0t^{-j}_n
\end{align*}
for any $n\geq 1$ and $0\leq j\leq n$.  With this definition at hand one
can easily see that $p_*$ is a morphism of cyclic modules.
\end{proof}

\begin{rem}
Let us see how each of the face morphisms of $\B{CM}^a_*(H,X)$ work:
First, let $n\geq 1$ and consider
\begin{align*}
t_n & (h^1\otimes\cdots\otimes h^n\otimes x)\\
 = & i_n\tau_np_n(h^1\otimes\cdots\otimes h^n\otimes x)\\
 = & i_n\tau_n\left(h^1_{(1)}\otimes\cdots\otimes h^n_{(1)}\otimes
            x_{(-1)}S^{-1}(h^1_{(3)}\cdots h^n_{(3)})
            \otimes h^1_{(2)}\cdots h^n_{(2)}x_{(0)}\right)\\
 = & i_n\left(x_{(-1)(1)}S^{-1}(h^1_{(3)(2)}\cdots h^n_{(3)(2)})
            \otimes h^1_{(1)}\otimes\cdots\otimes h^n_{(1)}
            \otimes x_{(-1)(2)}S^{-1}(h^1_{(3)(1)}\cdots h^n_{(3)(1)})
                    h^1_{(2)}\cdot h^n_{(2)}x_{(0)}\right)\\
 = & i_n\left(x_{(-1)}S^{-1}(h^1_{(2)}\cdots h^n_{(2)})
            \otimes h^1_{(1)}\otimes\cdots\otimes h^n_{(1)}\otimes
            x_{(0)}\right)\\
 = & \left(x_{(-1)}S^{-1}(h^1_{(2)}\cdots h^n_{(2)})
           \otimes h^1_{(1)}\otimes\cdots\otimes
           h^n_{(1)}x_{(0)}\right)
\end{align*}
Note that 
\begin{align*}
d_j(h^1\otimes\cdots\otimes h^n\otimes x)
 = & t_{n-1}^jd_0t_n^{-j}(h^1\otimes\cdots\otimes h^n\otimes x)\\
 = & i_{n-1}\tau_{n-1}^j\partial_0\tau_n^{-j}p_n
     (h^1\otimes\cdots\otimes h^n\otimes x)\\
 = & i_{n-1}\partial_j\left(
          h^1_{(1)}\otimes\cdots\otimes h^n_{(1)}
          \otimes x_{(-1)}S^{-1}(h^1_{(3)}\cdots\otimes h^n_{(3)})
          \otimes h^1_{(2)}\cdots h^n_{(2)}x_{(0)}\right)\\
 = &  \begin{cases}
      (\cdots\otimes h^{j+1}h^{j+2}\otimes\cdots\otimes x)
          & \text{ if } 0\leq j<n-1\\
      (h^1\otimes\cdots\otimes h^nx)
          & \text{ if } j=n-1\\
      \left(x_{(-1)}S^{-1}(h^1_{(2)}\cdots h^n_{(2)})h^1_{(1)}\otimes
                   h^2_{(1)}\otimes\cdots\otimes h^{n-1}_{(1)}\otimes
                   h^n_{(1)}x_{(0)}\right)
          & \text{ if } j=n
      \end{cases}
\end{align*}
\end{rem}

\begin{thm}\label{bar_isomorphism}
Let $X$ be a stable anti-Yetter-Drinfeld module/comodule and let
$BC^a_*(k,H,X)$ be the graded $k$--module $\{H^{\otimes n}\otimes
X\}_{n\geq 0}$ given with the cyclic structure
\begin{align}
\delta_j(h^1\otimes\cdots\otimes h^n\otimes x)
 = & \begin{cases}
     \epsilon(h^1)(h^2\otimes\cdots\otimes h^n\otimes x)
         & \text{ if } j=0\\
     (\cdots\otimes h^jh^{j+1}\otimes \cdots\otimes x)
         & \text{ if } 0<j<n\\
     (h^1\otimes\cdots\otimes h^{n-1}\otimes h^nx)
         & \text{ if } j=n
     \end{cases}\\
t_n(h^1\otimes\cdots\otimes h^n\otimes x)
 = & \left(x_{(-1)}S^{-1}(h^1_{(2)}\cdots h^n_{(2)})
           \otimes h^1_{(1)}\otimes\cdots\otimes
           h^n_{(1)}x_{(0)}\right)
\end{align}
Then there is an isomorphism of cyclic modules
$BC^a_*(k,H,X)\xra{t_*}\B{CM}^a_*(H,X)$.
\end{thm}

\begin{proof}
We need to check $t_*\delta_j=d_jt_*$ and $t_*^jt_*=t_*t_*^j$ for all
possible $j\in\B{Z}$.  The latter assertion is obvious.  For the former
observe that for $0<j\leq n$
\begin{align*}
d_j t_n = t_{n-1}^jd_0t_n^{-j+1}  = t_{n-1} d_{j-1} = t_{n-1}\delta_j
\end{align*}
and finally for $j=0$
\begin{align*}
t_{n-1}\delta_0(h^1\otimes\cdots\otimes h^n\otimes x)
 = & t_{n-1}\epsilon(h^1)(h^2\otimes\cdots\otimes h^n\otimes x)\\
 = & \epsilon(h^1)\left(x_{(-1)}S^{-1}(h^2_{(2)}\cdots h^n_{(2)})
       \otimes h^2_{(1)}\otimes\cdots\otimes h^{n-1}_{(2)}\otimes
        h^n_{(2)}x_{(0)}\right)\\
 = & d_0t_n(h^1\otimes\cdots\otimes h^n\otimes x)
\end{align*}
as we wanted to prove.
\end{proof}

\begin{rem}
By using Theorem~\ref{bar_isomorphism} one can conclude that, the
Hochschild complex of the cyclic module $\B{CM}^a_*(H,X)$ is isomorphic to
the bar complex $B^a_*(k,H,X)$.
\end{rem}

\begin{lem}\label{main_lemma}
Assume $Y$ is a right $H$--comodule and $X$ is a left $H$--comodule.
Endow $Y\otimes X$ with the following $H$--comodule structure:
\begin{align*}
\rho_R(y\otimes x)
 = (y_{(0)}\otimes x_{(0)})\otimes y_{(1)}S(x_{(-1)})
\end{align*}
for all $(y\otimes x)\in Y\otimes X$.  Then the cotensor product
\begin{align*}
Y\fotimes{H}X = \left\{(y\otimes x)\in Y\otimes X|\ 
    y_{(0)}\otimes y_{(1)}\otimes x = y\otimes x_{(-1)}\otimes x_{(0)}
    \right\}
\end{align*}
is isomorphic to
\begin{align*}
(Y\otimes X)^H=\left\{(y,x)\in Y\otimes X|\ 
  \rho_R(y\otimes x)=(y\otimes x\otimes\B{I})\right\}
\end{align*}
\end{lem}

\begin{proof}
First, let me show that $\rho_R$ is a genuine $H$--comodule structure:
\begin{align*}
(\rho_R\otimes id_H)\rho_R(y\otimes x)
 = & \left((y_{(0)(0)}\otimes x_{(0)(0)})\otimes y_{(0)(1)}S(x_{(0)(-1)})
           \otimes y_{(1)}S(x_{(-1)})\right)\\
 = & \left((y_{(0)}\otimes x_{(0)})\otimes y_{(1)}S(x_{(-1)})
            \otimes y_{(2)}S(x_{(-2)})\right)\\
 = & \left((y_{(0)}\otimes x_{(0)})\otimes y_{(1)(1)}S(x_{(-1)(2)})
            \otimes y_{(1)(2)}S(x_{(-1)(1)})\right)\\
 = & (id_{Y\otimes X}\otimes\Delta)\rho_R(y\otimes x)
\end{align*}
Now, consider $(y\otimes x)$ from $Y\fotimes{H} X$ and consider
\begin{align*}
\rho_R(y\otimes x) 
 = \left(y_{(0)}\otimes x_{(0)}\otimes y_{(1)}S(x_{(-1)})\right)
\end{align*}
But since $(y\otimes x)$ is form $Y\fotimes{H}X$, we have
$(y_{(0)}\otimes y_{(1)}\otimes x_{(-1)}\otimes x_{(0)}) =
(y_{(0)}\otimes y_{(1)}\otimes y_{(2)}\otimes x)$ which means
$\rho_R(y\otimes x) = (y\otimes x\otimes\B{I})$.  Conversely, if
$y\otimes x$ is from $(Y\otimes X)^H$ then
\begin{align*}
\left(y\otimes \B{I}x_{(-1)}\otimes x_{(0)}\right)
 = \left(y_{(0)}\otimes y_{(1)}S(x_{(-1)})x_{(0)(-1)}
        \otimes x_{(0)(0)}\right)
 = (y_{(0)}\otimes y_{(1)}\otimes x)
\end{align*}
which shows $Y\fotimes{H}X=(Y\otimes X)^H$.
\end{proof}

\begin{lem}\label{graded_coaction}
$\B{T}^a_*(H,X)$ is a $H$--comodule with the coaction defined as
\begin{align*}
\rho_R(h^0\otimes\cdots\otimes h^n\otimes x)
 =  (h^0_{(1)}\otimes\cdots\otimes h^n_{(1)}\otimes x_{(0)})
     \otimes h^0_{(2)}\cdots h^n_{(2)}S(x_{(-1)})
\end{align*}
$(h^0\otimes\cdots\otimes h^n\otimes x)$ from $\B{T}^a_*(H,X)$.
\end{lem}

\begin{thm}
Assume $X$ is a stable anti-Yetter-Drinfeld module.  Then, the morphism
of para-cyclic modules $\B{CM}^a_*(H,X)\xra{p_*}\B{T}^a_*(H,X)$ factors as
\begin{align*}
\B{CM}^a_*(H,X)\xra{p_*}\B{T}^a_*(H,X)^H\xra{id_*}\B{T}^a_*(H,X)
\end{align*}
Moreover, $p_*$ is an isomorphism of para-cyclic $H$--comodules.
\end{thm}

\begin{proof}
Let $z=h^1\cdots h^n$ and consider the expression
\begin{align*}
\rho_Rp_n & (h^0\otimes\cdots\otimes h^n\otimes x)\\
 = & \rho_R\left(h^1_{(1)}\otimes\cdots\otimes h^n_{(1)}\otimes
      x_{(-1)}S^{-1}(h^1_{(3)}\cdots h^n_{(3)})\otimes
      h^1_{(2)}\cdots h^n_{(2)}x_{(0)}\right)\\
 = & \rho_R\left(h^1_{(1)}\otimes\cdots\otimes h^n_{(1)}\otimes
      S(z_{(2)})z_{(3)}x_{(-1)}S^{-1}(z_{(5)})\otimes z_{(4)}x_{(0)}\right)\\
 = & \rho_R\left(h^1_{(1)}\otimes\cdots\otimes h^n_{(1)}\otimes
      S(z_{(2)})(z_{(3)}x)_{(-1)}\otimes (z_{(3)}x)_{(0)}\right)\\
 = & \left(h^1_{(1)(1)}\otimes\cdots\otimes h^n_{(1)(1)}\otimes
          \left(S(z_{(2)})(z_{(3)}x)_{(-1)}\right)_{(1)}
          \otimes (z_{(3)}x)_{(0)(0)}\right)\\
   & \hspace{1.5cm}\otimes\left(z_{(1)}S(z_{(2)})(z_{(3)}x)_{(-1)}
           \right)_{(2)} S\left((z_{(3)}x)_{(0)(-1)}\right)\\
 = & \left(h^1_{(1)(1)}\otimes\cdots\otimes h^n_{(1)(1)}\otimes
          S(z_{(2)(2)})(z_{(3)}x)_{(-1)(1)}
          \otimes (z_{(3)}x)_{(0)(0)}\right)\\
   & \hspace{1.5cm}\otimes 
          z_{(1)(2)}S(z_{(2)(1)})(z_{(3)}x)_{(-1)(2)}
          S\left((z_{(3)}x)_{(0)(-1)}\right)\\
 = & \left(h^1_{(1)}\otimes\cdots\otimes h^n_{(1)}\otimes
          S(z_{(2)})(z_{(3)}x)_{(-3)}
          \otimes (z_{(3)}x)_{(0)}\right)\otimes 
     (z_{(3)}x)_{(-2)}S\left((z_{(3)}x)_{(-1)}\right)\\
 = & \left(h^1_{(1)}\otimes\cdots\otimes h^n_{(1)}\otimes
          S(z_{(2)})(z_{(3)}x)_{(-1)}
          \otimes (z_{(3)}x)_{(0)}\right)\otimes\B{I}\\
 = & \left(h^1_{(1)}\otimes\cdots\otimes h^n_{(1)}\otimes
          S(z_{(2)})z_{(3)(1)}x_{(-1)}S^{-1}(z_{(3)(3)})
          \otimes z_{(3)(2)}x_{(0)}\right)\otimes\B{I}\\
 = & \left(h^1_{(1)}\otimes\cdots\otimes h^n_{(1)}\otimes
          x_{(-1)}S^{-1}(z_{(3)})
          \otimes z_{(2)}x_{(0)}\right)\otimes\B{I}\\
 = & \left(h^1_{(1)}\otimes\cdots\otimes h^n_{(1)}\otimes
          x_{(-1)}S^{-1}(h^1_{(3)}\cdots h^n_{(3)})\otimes
          h^1_{(2)}\cdots h^n_{(2)}x\right)\otimes\B{I}\\
 = & p_n(h^0\otimes\cdots\otimes h^n\otimes x)\otimes\B{I}
\end{align*}
which means $p_*$ factors as
$\B{CM}^a_*(H,X)\xra{p_*}\B{T}^a_*(H,X)^H\xra{id_*}\B{T}^a_*(H,X)$ as we
wanted to prove.

With the help Lemma~\ref{main_lemma}, we can say graded submodule
$\B{T}^a_*(H,X)^H$ consists of elements of the form
$(h^0\otimes\cdots\otimes h^n\otimes x)$ such that
\begin{align*}
(h^0_{(1)}\otimes\cdots\otimes h^n_{(1)}\otimes h^1_{(2)}\cdots
  h^n_{(2)}\otimes x)
 = (h^0\otimes\cdots\otimes h^n\otimes x_{(-1)}\otimes x_{(0)})
\end{align*}
Take $(h^1\otimes\cdots\otimes h^n\otimes x)$ from $\B{T}(H,X)^H$ and
consider
\begin{align*}
p_ni_n
   & (h^1\otimes\cdots\otimes h^{n+1}\otimes x)\\
 = & p_n\left(h^1\otimes\cdots\otimes h^n\otimes h^{n+1}x\right)\\
 = & \left(h^1_{(1)}\otimes\cdots\otimes h^n_{(1)}\otimes
           h^{n+1}_{(1)}x_{(-1)}S^{-1}(h^{n+1}_{(3)})
           S^{-1}(h^1_{(3)}\cdots h^n_{(3)})\otimes
           h^1_{(2)}\cdots h^n_{(2)}h^{n+1}_{(2)}x_{(0)}\right)\\
 = & \left(h^1_{(1)}\otimes\cdots\otimes h^n_{(1)}\otimes
           h^{n+1}_{(1)}x_{(-1)}
           S^{-1}(h^1_{(3)}\cdots h^n_{(3)}h^{n+1}_{(3)})\otimes
           h^1_{(2)}\cdots h^n_{(2)}h^{n+2}_{(2)}x_{(0)}\right)\\
 = & \left(h^1_{(1)}\otimes\cdots\otimes h^n_{(1)}\otimes
           h^{n+1}_{(1)}h^1_{(4)}\cdots h^{n+1}_{(4)}
           S^{-1}(h^1_{(3)}\cdots h^{n+1}_{(3)})\otimes
           h^1_{(2)}\cdots h^{n+1}_{(2)}x\right)\\
 = & \left(h^1_{(1)}\otimes\cdots\otimes h^n_{(1)}\otimes
           h^{n+1}_{(1)}\otimes
           h^1_{(2)}\cdots h^{n+1}_{(2)}x\right)\\
 = & (h^1\otimes\cdots\otimes h^{n+1}\otimes x_{(-1)}x_{(0)})\\
 = & (h^1\otimes\cdots\otimes h^{n+1}\otimes x)
\end{align*}
which means $p_*$ is an epimorphism.  We already had $id_*=i_*p_*$ from
Remark~\ref{monomorphism} which makes $p_*$ a monomorphism too.  Thus
$p_*$ is an isomorphism.
\end{proof}

\begin{thm}
Assume $X$ is a stable anti-Yetter-Drinfeld module/comodule.  Then
$\B{CM}^a_*(H,X)$ is a cyclic $k$--module and
$\B{CM}^a_*(H,X)\xra{p_*}\B{T}^a_*(H,X)^H$ is an isomorphism of cyclic
$k$--modules.
\end{thm}

\begin{proof}
Observe that, if we assume $X$ is a $0$--stable module/comodule and
$(h^0\otimes\cdots\otimes h^n\otimes x)$ is in $\B{T}^a_n(H,X)^H$, then
\begin{align*}
\tau_n^{n+1}(h^0\otimes\cdots\otimes h^n\otimes x)
 = & (h^0_{(1)}\otimes\cdots\otimes h^n_{(1)}\otimes 
      h^0_{(2)}\cdots h^n_{(2)}x)\\
 = & (h^0\otimes\cdots\otimes h^n\otimes x_{(-1)}x_{(0)})\\
 = & (h^0\otimes\cdots\otimes h^n\otimes x)
\end{align*}
Therefore, if $X$ is a stable anti-Yetter-Drinfeld module/comodule then
the para-cyclic $k$--submodule of $H$--coaction invariants
$\B{T}^a_*(H,X)^H$ of $\B{T}^a_*(H,X)$ is a cyclic module and
$\B{CM}^a_*(H,X)\xra{p_*}\B{T}^a_*(H,X)^H$ is an isomorphism of cyclic
modules.
\end{proof}

\section{Bialgebra cyclic homology}\label{NewTheory}

\begin{defn}
For any $j\in\B{Z}$, define a degree 1 morphism $[\rho_R,\tau_*^j]$ on
$\B{T}^a_*(H,X)$ as
\begin{align*}
[\rho_R,\tau_*^j] = & 
  \rho_R\tau_n^j - (\tau_n^j\otimes\B{I})\rho_R 
\end{align*}
and let $\B{PCM}^a_*(H,X):= \bigcap_{j\in\B{Z}}ker([\rho_R,\tau_*^j])$.
\end{defn}

\begin{thm}\label{cyclic}
$\B{PCM}^a_*(H,X)$ is a para-cyclic $H$--comodule whenever $H$ is a Hopf
algebra and $X$ is an anti-Yetter-Drinfeld module.  Moreover,
$\B{CM}^a_*(H,X)$ is isomorphic to the cyclic $H$--subcomodule
$\B{PCM}^a_*(H,X)^H$.
\end{thm}

\begin{proof}
Note that for any $p_*({\bf h}\otimes x)$ in the image of $p_*$, we have
\begin{align*}
[\rho_R,\tau_*^j]p_*({\bf h}\otimes x)
 = & \rho_R\tau_*^jp_*({\bf h}\otimes x)
      - (\tau_*^j\otimes id_H)\rho_Rp_*({\bf h}\otimes x)\\
 = & \rho_Rp_*t_*^j({\bf h}\otimes x)
      - \tau_*^jp_*({\bf h}\otimes x)\otimes\B{I}\\
 = & p_*t_*^j({\bf h}\otimes x)\otimes\B{I}
      - \tau_*^jp_*({\bf h}\otimes x)\otimes\B{I}
 = 0 
\end{align*}
which means $\B{PCM}^a_*(H,X)\supseteq im(p_*)=\B{T}^a_*(H,X)^H$.

The graded $k$--submodule $\B{PCM}^a_*(H,X)$ is stable under the actions
of $\tau_*^i$ for any $i\in\B{Z}$ since for any $({\bf h}\otimes x)$
from $\B{PCM}^a_*(H,X)$ we have
\begin{align*}
[\rho_R,\tau_*^j]\tau_*^i({\bf h}\otimes x)
 = -\tau_*^j[\rho_R,\tau_*^i]({\bf h}\otimes x)
   +[\rho_R,\tau_*^{j+i}]({\bf h}\otimes x)
 = 0 
\end{align*}
We need to show that $\B{PCM}^a_*(H,X)$ is actually a para-cyclic
submodule of $\B{T}^a_*(H,X)$.  In order to prove this, we need to prove
that $\B{PCM}_n(H,X)$ is stable under the action of $\partial_0$ for any
$n\geq 0$.

First, observe that for $0\leq j\leq n-1$
\begin{align*}
\rho_R\partial_j(h^0\otimes\cdots\otimes h^n\otimes x)
 = & \rho_R(\cdots\otimes h^jh^{j+1}\otimes\cdots\otimes x)\\
 = & \left(\cdots\otimes h^j_{(1)}h^{j+1}_{(1)}
           \otimes\cdots\otimes x_{(0)}\right)\otimes
     h^0_{(2)}\cdots h^n_{(2)}S(x_{(-1)})\\
 = & (\partial_j\otimes id_H)\rho_R(h^0\otimes\cdots\otimes 
      h^n\otimes x)
\end{align*}
and for $j=n$, assuming $(h^0\otimes\cdots\otimes h^n\otimes x)$ is in
$\B{PCM}^a_*(H,X)$ we obtain
\begin{align*}
\rho_R\partial_n(h^0\otimes\cdots\otimes h^n\otimes x)
 = & \rho_R\partial_0\tau_n(h^0\otimes\cdots\otimes h^n\otimes x)\\
 = & (\partial_0\otimes id_H)\rho_R\tau_n(h^0\otimes\cdots\otimes
      h^n\otimes x)
\end{align*}
However, since $[\rho_R,\tau_*](h^0\otimes\cdots\otimes h^n\otimes
x)=0$, we have
\begin{align*}
\rho_R\partial_n(h^0\otimes\cdots\otimes h^n\otimes x)
 = & (\partial_0\otimes id_H)(\tau_n\otimes id_H)\rho_R
        (h^0\otimes\cdots\otimes h^n\otimes x)\\
 = & (\partial_n\otimes id_H)\rho_R        
        (h^0\otimes\cdots\otimes h^n\otimes x)
\end{align*}
Notice also that, the definition dictates that
\begin{align*}
\tau_{n-1}^j \partial_0 
  = \partial_{j\text{ mod } n+1}\tau_n^j 
\end{align*}
for any $j\in\B{Z}$ after observing the fact that $\partial_0 =
\tau_{n-1}^n\partial_0\tau_n^{-n-1}$.

Now, consider $[\rho_R,\tau_*^j ]\partial_0$ restricted to
$\B{PCM}^a_*(H,X)$ 
\begin{align*}
[\rho_R,\tau_*^j ]\partial_0
 = & \rho_R\tau_*^j \partial_0 
     - (\tau_*^j\otimes id_H)\rho_R\partial_0\\
 = & \rho_R\partial_j\tau_*^j  
     - (\tau_*^j\partial_0\otimes id_H)\rho_R\\
 = & (\partial_j\otimes id_H)\rho_R\tau_*^j  
     - (\partial_j\tau_*^j\otimes id_H)\rho_R\\
 = & (\partial_j\otimes id_H)[\rho_R,\tau_*^j ]
\end{align*}
which is uniformly zero, in other words $\B{PCM}^a_*(H,X)$ is stable under
the action of $\partial_0$ as we wanted to show.  This finishes the the
proof that $\B{PCM}^a_*(H,X)$ is a para-cyclic submodule of
$\B{T}^a_*(H,X)$.

Now, let me show that $\B{PCM}^a_*(H,X)$ is a graded
$H$--subcomodule of $\B{T}^a_*(H,X)$.  For this end, consider 
\begin{align*}
([\rho_R,\tau_*^j]\otimes id_H)\rho_R
 = & (\rho_R\tau_*^j\otimes id_H)\rho_R
      - (\tau_*^j\otimes id_{H\otimes H})(\rho_R\otimes id_H)\rho_R\\
 = & (\rho_R\tau_*^j\otimes id_H)\rho_R 
      - (\tau_*^j\otimes id_{H\otimes H})(id_*\otimes\Delta)\rho_R\\
 = & (\rho_R\tau_*^j\otimes id_H)\rho_R 
      - (\tau_*^j\otimes\Delta)\rho_R
\end{align*}
Restricted to $\B{PCM}^a_*(H,X)$, $(\tau_*^j\otimes
id_H)\rho_R=\rho_R\tau_*$.  Thus
\begin{align*}
([\rho_R,\tau_*^j]\otimes id_H)\rho_R
 = & (\rho_R\tau_*^j\otimes id_H)\rho_R 
      - (id_*\otimes\Delta_H)\rho_R\tau_*^j\\
 = & (\rho_R\tau_*^j\otimes id_H)\rho_R 
      - (\rho_R\otimes id_H)\rho_R\tau_*^j\\
 = & (\rho_R\otimes id_H)[\rho_R,\tau_*]
\end{align*}
is uniformly zero.  Therefore, $\rho_R$ sends $\B{PCM}^a_*(H,X)$ to
$H\otimes \B{PCM}^a_*(H,X)$, i.e. it is a $H$--comodule.

Now that we showed $\B{PCM}^a_*(H,X)$ is a para-cyclic module and
$H$--comodule, let me merge these two structures and show that it is a
para-cyclic $H$--comodule: For this end, we must show that the
$H$--coaction and the action of the cyclic groups and the primary face
maps $\partial_0$ commute.  We already proved that the $H$--coaction and
$\partial_0$ commute on $\B{PCM}^a_*(H,X)$.  Moreover, the $H$--coaction
and the action of cyclic groups commute by design on $\B{PCM}^a_*(H,X)$.

Finally observe that since $\B{T}^a_*(H,X)^H\subseteq
\B{PCM}^a_*(H,X)\subseteq \B{T}^a_*(H,X)$ is a chain of monomorphisms of
graded $H$--comodules, we have
\begin{align*}
\B{T}^a_*(H,X)^H = \B{PCM}^a_*(H,X)^H \cong \B{CM}^a_*(H,X)
\end{align*}
\end{proof}

\begin{thm}
Let $X$ be a stable anti-Yetter-Drinfeld module/comodule and let $H$ be
a commutative Hopf algebra.  Then
$\B{PCM}^a_*(H,X)=\B{T}^a_*(H,X)$.
\end{thm}

\begin{proof}
Assume $(h^0\otimes\cdots\otimes h^n\otimes x)$ is an arbitrary element
of $\B{T}^a_*(H,X)$.  Consider
\begin{align*}
\rho_R\tau_n^j 
   & (h^0\otimes\cdots\otimes h^n\otimes x)\\
 = & \rho_R\left(h^j\otimes\cdots\otimes h^n\otimes h^0_{(1)}
      \otimes\cdots\otimes
      h^{j-1}_{(1)}\otimes h^0_{(2)}\cdots h^{j-1}_{(2)}x\right)\\
 = & \left(h^j_{(1)}\otimes\cdots\otimes h^n_{(1)}\otimes h^0_{(1)(1)}
      \otimes\cdots\otimes
      h^{j-1}_{(1)(1)}\otimes h^0_{(2)(2)}\cdots
      h^{j-1}_{(2)(2)} x_{(0)}\right)\\
   & \hspace{1cm}\otimes
      h^j_{(2)}\cdots h^n_{(2)}h^0_{(1)(2)}\cdots h^{j-1}_{(1)(2)}
      h^0_{(2)(3)}\cdots h^{j-1}_{(2)(3)}S(x_{(-1)})
      S(h^0_{(2)(1)}\cdots h^{j-1}_{(2)(1)})\\
 = & \left(h^j_{(1)}\otimes\cdots\otimes h^n_{(1)}\otimes h^0_{(1)}
      \otimes\cdots\otimes
      h^{j-1}_{(1)}\otimes h^0_{(4)}\cdots
      h^{j-1}_{(4)} x_{(0)}\right)\\
   & \hspace{1cm}\otimes
      h^0_{(5)}\cdots h^{j-1}_{(5)}h^j_{(2)}\cdots 
      h^n_{(2)}h^0_{(2)}\cdots h^{j-1}_{(2)}
      S(h^0_{(3)}\cdots h^{j-1}_{(3)})S(x_{(-1)})\\
 = & \left(h^j_{(1)}\otimes\cdots\otimes h^n_{(1)}\otimes h^0_{(1)}
      \otimes\cdots\otimes h^{j-1}_{(1)}\otimes h^0_{(2)}\cdots
      h^{j-1}_{(2)} x_{(0)}\right)
      \otimes h^0_{(3)}\cdots h^{j-1}_{(3)}
              h^j_{(2)}\cdots h^n_{(2)}S(x_{(-1)})
\end{align*}
Then
\begin{align*}
(\tau_n^{-1}\otimes id_H)
 \rho_R\tau_n^j(h^0\otimes\cdots\otimes h^n\otimes x)
 = & (h^0_{(1)}\otimes\cdots\otimes h^n_{(1)}\otimes x_{(0)})
     \otimes h^0_{(2)}\cdots h^n_{(2)}S(x_{(-1)})\\
 = & \rho_R(h^0\otimes\cdots\otimes h^n\otimes x)
\end{align*}
which means $[\rho_R,\tau_*]\equiv 0$ uniformly on $\B{T}^a_*(H,X)$.
\end{proof}

\begin{cor}
Let $H$ and $X$ be as before.  Then $\B{PCM}^a_*(H,X)\cong H\otimes
BC^a_*(k,H,X)$.
\end{cor}

\begin{proof}
Since $H$ is commutative $ad(H)$ is a trivial $H$--module.  Therefore
\begin{align*}
\B{T}^a_*(H,X)
\cong BC^a_*(ad(H),H,X) \cong H\otimes BC^a_*(k,H,X)
\end{align*}
as we wanted to prove.
\end{proof}

\begin{rem}
Any left $H$--comodule $M$ is isomorphic to a right $H$--comodule
$M^{op}$ as follows:  Define
\begin{align*}
\rho_R^M(m)= m_{(0)}\otimes S(m_{(-1)})
\end{align*}
Similarly any right $H$--comodule $N$ is isomorphic to a left
$H$--comodule $N^{op}$ via 
\begin{align*}
\rho_R^N(n) = S^{-1}(n_{(1)})\otimes n_{(0)}
\end{align*}
One can immediately see that $(M^{op})^{op}=M$ and $(N^{op})^{op}=N$ for
any comodules $M$ and $N$.
\end{rem}

\begin{rem}\label{redefined}
Assume $X$ is a right $H$--comodule and assume we used
$\B{T}^a_*(H,X^{op})$ above.  Now, instead of the left $H$--comodule
structure $\rho_R$ on $\B{T}^a_*(H,X^{op})$ we defined above, we can use
\begin{align*}
\rho_R(h^0\otimes\cdots\otimes h^n\otimes x)
 = & \left(h^0_{(1)}\otimes\cdots\otimes h^n_{(1)}\otimes x_{(0)}
         \otimes h^0_{(2)}\cdots h^n_{(2)}) S(x_{(-1)})\right)\\
 = & \left(h^0_{(1)}\otimes\cdots\otimes h^n_{(1)}\otimes x_{(0)}
         \otimes h^0_{(2)}\cdots h^n_{(2)}x_{(1)}\right)
\end{align*}
The advantage of the these re-writings of the coaction is that, even 
when $B$ is just a bialgebra,  for a $B$--comodule algebra $Y$ and 
$0$--stable $B$--module/comodule $X$, now we can define $\B{PCM}^a_*(Y,X)$
and $\B{CM}^a_*(Y,X)$.
\end{rem}

\begin{thm}
Let $Y$ be a right $B$--comodule algebra and $X$ be a left $B$--module
and a right $B$--comodule.  Assume also that $X$ satisfies the property
that $x_{(1)}x_{(0)}=x$ for all $x\in X$.  Then $\B{T}^a_*(Y,X)$ is a
simplicial module which is short of being a para-cyclic modules since
$\tau_*$ may not be invertible.  In the case $B$ is a Hopf algebra, 
$X^{op}$ is $1$--stable, $\B{T}^a_*(Y,X^{op})$ is a para-cyclic module and 
the submodule $\B{PCM}^a_*(Y,X^{op})$ is a para-cyclic $B$--comodule.  
If we define $\B{CM}^a_*(Y,X)$ as the graded submodule of 
$\B{PCM}^a_*(Y,X^{op})$ containing elements of the form
$(y^0\otimes\cdots\otimes y^n\otimes x)\in
\B{PCM}^a_*(H,X)$ such that
\begin{align}\label{cotensor}
\rho_R(y^0\otimes\cdots\otimes y^n\otimes x) 
 = (y^0_{(0)}\otimes\cdots\otimes y^n_{(0)}\otimes x_{(0)})
    \otimes y^0_{(1)}\cdots y^n_{(1)}x_{(1)}
 = (y^0\otimes\cdots\otimes y^n\otimes x)\otimes\B{I}
\end{align}
then $\B{CM}^a_*(Y,X)$ is a cyclic module, regardless of $B$ being a Hopf
algebra.  In the case of $B$ is a Hopf algebra $\B{CM}^a_*(Y,X)$ is the
same as $\B{PCM}^a_*(Y,X^{op})^B$.
\end{thm}

\begin{proof}
The graded $k$--module $\B{T}^a_*(Y,X)$ is the collection $\{Y^{\otimes
n+1}\otimes X\}_{n\geq 0}$.  The simplicial structure is given by the
structure maps defined in Equation~\ref{simplicial_definition} as
\begin{align}
\partial_j(y^0\otimes\cdots\otimes y^n\otimes x)
 = \begin{cases}
     (\cdots\otimes y^jy^{j+1}\otimes\cdots\otimes x)
        & \text{ if } 0\leq j< n\\
     \left(y^n_{(0)}y^0\otimes\cdots\otimes h^{n-1}\otimes 
           y^n_{(1)}x\right)
        & \text{ if } j=n
     \end{cases}
\end{align}
The cyclic maps are defined in Equation~\ref{cyclic_definition} as
\begin{align}
\tau_n(y^0\otimes\cdots\otimes y^n\otimes x)
 = & \left(y^n_{(0)}\otimes y^0\otimes\cdots\otimes y^{n-1}\otimes 
       y^n_{(1)}x\right)
\end{align}
Since $B$ is just a bialgebra, we don't have an antipode.  Thus the
cyclic maps may not be invertible.  This means, $\B{T}^a_*(Y,X)$ is almost
a para-cyclic $k$--module, since all other identities are satisfied.
Regardless of $B$ being a Hopf algebra, one can still define
$\B{PCM}^a_*(Y,X)$ by using $\rho_R$ defined in Remark~\ref{redefined}
as
\begin{align}
\rho_R(y^0\otimes\cdots\otimes y^n\otimes x)
 = & (y^0_{(0)}\otimes\cdots\otimes y^n_{(0)}\otimes x_{(0)})\otimes 
     y^0_{(1)}\cdots y^n_{(1)}x_{(1)}
\end{align}
$\B{PCM}^a_*(Y,X)$ is still a para-cyclic $B$--module thanks to the
identity
\begin{align}
\tau_{n-1}^j\partial_0 = \partial_{j\text{ mod } n}\tau_n^j
\end{align}
for all $j\geq 0$.  $\B{CM}^a_*(Y,X)$ is always cyclic $B$--module since
$\rho_R(y^0\otimes\cdots\otimes y^n\otimes x) = (y^0\otimes\cdots\otimes
y^n\otimes x)\otimes\B{I}$ which implies
\begin{align*}
t_n^{n+1}(y^0\otimes\cdots\otimes y^n\otimes x)
 = & (y^0_{(0)}\otimes\cdots\otimes y^n_{(0)}\otimes 
      y^0_{(1)}\cdots y^n_{(1)}x)\\
 = & (y^0_{(0)}\otimes\cdots\otimes y^n_{(0)}\otimes 
      y^0_{(1)}\cdots y^n_{(1)}x_{(1)}x_{(0)})\\
 = & (y^0\otimes\cdots\otimes y^n\otimes x)
\end{align*}
for any $(y^0\otimes\cdots\otimes y^n\otimes x)\in \B{PCM}^a_*(Y,X)$

In the case of $B$ is a Hopf algebra, we use $X^{op}$.  Then the
condition given in Equation~\ref{cotensor} can be written as
\begin{align}
(y^0_{(0)}\otimes\cdots\otimes y^n_{(0)}\otimes 
 y^0_{(1)}\cdots y^n_{(1)}\otimes x)
 = (y^0\otimes\cdots\otimes y^n\otimes S^{-1}(x_{(1)})\otimes x_{(0)})
\end{align}
using Lemma~\ref{main_lemma} one can see that $\B{CM}^a_*(Y,X)\cong
\B{PCM}^a_*(Y,X^{op})^B$
\end{proof}

\section{Cyclic Duality}\label{Duality}

Assume $H$ is a Hopf algebra with a bijective antipode and $X$ is a
stable $H$--module/comodule.

\begin{rem}
In this section we need the full para-cyclic structure on
$\B{T}^a_*(H,X)$.  This means, we must provide degeneracy morphisms.
Define
\begin{align}
\sigma_0(h^0\otimes\cdots\otimes h^n\otimes x) 
 = & (h^0\otimes\B{I}\otimes h^1\otimes\cdots\otimes h^n\otimes x)\\
\sigma_j = & \tau_{n+1}^j\sigma_0\tau_n^{-j}
\end{align}
for all $0\leq j\leq n$ and for all $(h^0\otimes\cdots\otimes h^n\otimes
x)$ from $\B{T}^a_*(H,X)$.  We leave checking the cocyclic identities to
the reader.
\end{rem}

\begin{rem}
Recall from \cite{Kaygun:BialgebraCyclic} that, we have defined
$\B{T}^c_*(H,X)$ as the graded module $\{H^{\otimes n+1}\otimes
X\}_{n\geq 0}$ with the following cocyclic structure:
\begin{align*}
\partial^c_0(h^0\otimes\cdots\otimes h^n\otimes x)
 = & (h^0_{(1)}\otimes h^0_{(2)}\otimes h^1\otimes\cdots\otimes
      h^n\otimes x)\\
\tau_{c,n}(h^0\otimes\cdots\otimes h^n\otimes x)
 = & \left(S^{-1}(x_{(-1)})h^n\otimes h^0\otimes\cdots\otimes h^{n-1}
           \otimes x_{(0)}\right)\\
\partial^c_j = & \tau_{n+1}^j\partial_0\tau_n^{-j}
\end{align*}
for all $0\leq j\leq n+1$.  Now, to that add the codegeneracy maps which
are defined as
\begin{align*}
\sigma^c_0(h^0\otimes\cdots\otimes h^n\otimes x)
 = & \epsilon(h^1)(h^0\otimes h^2\otimes\cdots\otimes h^n\otimes x)\\
\sigma^c_j 
 = & \tau_{c,n-1}^j\sigma^c_0\tau_{c,n}^{-j}
\end{align*}
for all $0\leq j\leq n$ and $(h^0\otimes\cdots\otimes h^n\otimes x)$
from $\B{T}^c_*(H,X)$.  We leave checking the cocyclic identities to the
reader.
\end{rem}

\begin{lem}
Let $A_*=\{A_n\}_{n\geq 0}$ be a para-cocyclic $k$--module with
structure morphisms
\begin{align*}
A_{n-1} & \xra{\partial_j} A_n & 
A_{n+1} & \xra{\sigma_j}   A_n &
A_n     & \xra{\tau_n}     A_n 
\end{align*} 
for $0\leq j\leq n$.  Then the graded $k$--module $A_*^\vee =
\{A_n\}_{n\geq 0}$ with the following structure morphisms
\begin{align*}
A_n & \xra{\partial^\vee_0     = \sigma_{n-1}\tau_n} A_{n-1} & 
A_n & \xra{\partial^\vee_{i+1} = \sigma_j}           A_{n-1} & 
A_n & \xra{\sigma^\vee_j       = \partial_j}         A_{n+1} & 
A_n & \xra{\tau_{n,\vee}       = \tau_n^{-1}}        A_n 
\end{align*}
for all $0\leq j\leq n$ and $0\leq i\leq n-1$ is a para-cyclic
$k$--module.
\end{lem}

\begin{lem}\label{duality_lemma}
Define a morphism of graded modules
$\B{T}^c_*(H,X)\xra{\beta_*}\B{T}^a_*(H,X)$ by
\begin{align}
\beta_n(h^0\otimes\cdots\otimes h^n\otimes x)\otimes h
 = & \left(S(h^n_{(3)})x_{(-1)}h^0_{(1)}\otimes
           S(h^0_{(2)})h^1_{(1)}\otimes\cdots\otimes
           S(h^{n-1}_{(2)})h^n_{(1)}\otimes S(h^n_{(2)})x_{(0)}\right)
\end{align}
for all $(h^0\otimes\cdots\otimes h^n\otimes x)$ from $\B{T}^a_*(H,X)$.
Then one has
\begin{align}
\partial_1 \beta_n = & \beta_{n-1}\sigma^c_0    & 
\sigma_0   \beta_n = & \beta_{n+1}\partial^c_0  & 
\tau_n^{-1}\beta_n = & \beta_n    \tau_{c,n}  
\end{align}
for all $n\geq 0$.
\end{lem}

\begin{proof}
Proof is going to be by direct calculation.  For $n\geq 1$ consider
\begin{align*}
\partial_1\beta_n 
   & (h^0\otimes\cdots\otimes h^n\otimes x)\\
 = & \partial_1\left(S(h^n_{(3)})x_{(-1)}h^0_{(1)}\otimes
           S(h^0_{(2)})h^1_{(1)}\otimes\cdots\otimes
           S(h^{n-1}_{(2)})h^n_{(1)}\otimes S(h^n_{(2)})x_{(0)}\right)\\
 = & \begin{cases}
     \left(S(h^0_{(2)(2)})h^1_{(1)(1)}S(h^1_{(3)})x_{(-1)}h^0_{(1)}
           \otimes S(h^0_{(2)(1)})h^1_{(1)(2)}S(h^1_{(2)})x_{(0)}\right)
           & \text{ if } n=1\\
     \epsilon(h^1)\left(S(h^n_{(3)})x_{(-1)}h^0_{(1)}\otimes
           S(h^0_{(2)})h^2_{(1)}\otimes\cdots\otimes
           S(h^{n-1}_{(2)})h^n_{(1)}\otimes S(h^n_{(2)})x_{(0)}\right)
           & \text{ if } n\geq 2
     \end{cases}\\
 = & \beta_{n-1}\sigma^c_0(h^0\otimes\cdots\otimes h^n\otimes x)
\end{align*}
which proves the first identity.  Now let $n\geq 0$ and take
$(h^0\otimes\cdots\otimes h^n\otimes x)$ from $\B{T}^a_n(H,X)$ and
consider
\begin{align*}
\sigma_0\beta_n 
   & (h^0\otimes\cdots\otimes h^n\otimes x)\\
 = & \sigma_0\left(S(h^n_{(3)})x_{(-1)}h^0_{(1)}\otimes
           S(h^0_{(2)})h^1_{(1)}\otimes\cdots\otimes
           S(h^{n-1}_{(2)})h^n_{(1)}\otimes S(h^n_{(2)})x_{(0)}\right)\\
 = & \left(S(h^n_{(3)})x_{(-1)}h^0_{(1)}\otimes\B{I}\otimes
           S(h^0_{(2)})h^1_{(1)}\otimes\cdots\otimes
           S(h^{n-1}_{(2)})h^n_{(1)}\otimes S(h^n_{(2)})x_{(0)}\right)\\
 = & \left(S(h^n_{(3)})x_{(-1)}h^0_{(1)(1)}\otimes 
           S(h^0_{(1)(2)})h^0_{(2)(1)}\otimes
           S(h^0_{(2)(2)})h^1_{(1)}\otimes\cdots\otimes
           S(h^{n-1}_{(2)})h^n_{(1)}\otimes S(h^n_{(2)})x_{(0)}\right)\\
 = & \beta_{n+1}\partial^c_0(h^0\otimes\cdots\otimes h^n\otimes x)
\end{align*}
which proves the second identity.  Finally for $n\geq 1$
\begin{align*}
\tau_n^{-1}\beta_n
   & (h^0\otimes\cdots\otimes h^n\otimes x)\\
 = & \tau_n^{-1}\left(S(h^n_{(3)})x_{(-1)}h^0_{(1)}\otimes
           S(h^0_{(2)})h^1_{(1)}\otimes\cdots\otimes
           S(h^{n-1}_{(2)})h^n_{(1)}\otimes S(h^n_{(2)})x_{(0)}\right)\\
 = & \left(S(h^0_{(2)})h^1_{(1)}\otimes\cdots\otimes
           S(h^{n-1}_{(2)})h^n_{(1)}\otimes 
           S(h^n_{(3)(2)})x_{(-1)(1)}h^0_{(1)(1)}\right.\\
   & \hspace{1.5cm}\left.\otimes
           S(h^0_{(1)(2)})S(x_{(-1)(2)})S^2(h^n_{(3)(1)})
           S(h^n_{(2)})x_{(0)}\right)\\
 = & \left(S(h^0_{(3)})h^1_{(1)}\otimes\cdots\otimes
           S(h^{n-1}_{(2)})h^n_{(1)}\otimes 
           S(h^n_{(2)})x_{(-1)}h^0_{(1)}\otimes
           S(h^0_{(2)})x_{(0)}\right)\\
 = & \left(S(x_{(-1)(3)}h^0_{(3)})x_{(0)(-1)}h^1_{(1)} 
           \otimes\cdots\otimes
           S(h^{n-1}_{(2)})h^n_{(1)}\otimes 
           S(h^n_{(2)})x_{(-1)(1)}h^0_{(1)}\otimes
           S(x_{(-1)(2)}h^0_{(2)})x_{(0)(0)}\right)\\
 = & \beta_n\left(h^1\otimes\cdots\otimes h^n\otimes
           x_{(-1)}h^0\otimes x_{(0)}\right)\\
 = & \beta_n\tau_{c,n}(h^0\otimes\cdots\otimes h^n\otimes x)
\end{align*}
which proves the third identity.
\end{proof}

\begin{thm}
There is a morphism of para-cyclic modules
$\B{T}^c_*(H,X)^\vee\xra{\beta_*}\B{T}^a_*(H,X)$.
\end{thm}
\begin{proof}
By using Lemma~\ref{duality_lemma} one can see that
\begin{align*}
\partial_{i+1}\beta_n 
 = & \beta_{n-1}\sigma^c_j = \beta_{n-1}\partial^{c,\vee}_{j+1}\\
\sigma_j\beta_n
 = & \beta_{n+1}\partial^c_j = \beta_{n+1}\sigma^{c,\vee}_j\\
\tau_n \beta_n
 = & \beta_n\tau_{c,n}^{-1} = \beta_n\tau_{c,n,\vee}
\end{align*}
for any $0\leq i\leq n$ and $0\leq j\leq n+1$.  This also implies
\begin{align*}
\beta_n\partial^{c,\vee}_0 
 = & \beta_n\sigma^c_{n-1}\tau_{c,n}
 = \partial_n\tau_n^{-1}\beta_n = \partial_0\beta_n
\end{align*}
for any $n\geq 1$.
\end{proof}

\begin{defn}
Define a morphism of graded modules
$\B{T}^a_*(H,X)\xra{\alpha_*}\B{T}^c_*(H,X)^\vee$ by letting
\begin{align*}
\alpha_n(h^0\otimes\cdots\otimes h^n\otimes x)
 = & \left(h^0_{(1)}\otimes h^0_{(2)}h^1_{(1)}
           \otimes\cdots\otimes 
           h^0_{(n+1)}h^1_{(n)}\cdots h^n_{(1)}\otimes
           h^0_{(n+2)}h^1_{(n+1)}\cdots h^n_{(2)}x\right)
\end{align*}
for all $(h^0\otimes\cdots\otimes h^n\otimes x)$ from $\B{T}^c_n(B,X)$
for an arbitrary $n\geq 0$. 
\end{defn}

\begin{lem}\label{graded_epimorphism}
$\beta_*\alpha_* = id_*$ restricted to $\B{CM}^a_*(H,X)$.
\end{lem}

\begin{proof}
For $n=0$, one has
\begin{align*}
\alpha_0 (h^0\otimes x) = & (h^0_{(1)}\otimes h^0_{(2)}x)
\end{align*}
But, recall that we took $(h^0\otimes x)$ from $\B{CM}^a_*(H,X)$.  Then
\begin{align*}
(h^0_{(1)}\otimes h^0_{(2)}\otimes x)
 = (h^0\otimes x_{(-1)}\otimes x_{(0)})
\end{align*}
which implies
\begin{align*}
\alpha_0 (h^0\otimes x)
 = & (h^0\otimes x_{(-1)}x_{(0)})
 = (h^0\otimes x) 
\end{align*}
Then
\begin{align*}
\beta_0\alpha_0 (h^0\otimes x)
 = \beta_0(h^0\otimes x)
 = \left(S(h^0_{(3)})x_{(-1)}h^0_{(1)}\otimes S(h^0_{(2)})x_{(0)}\right)
\end{align*}
Again by using the fact that $(h^0\otimes x)$ is from $\B{CM}^a_*(H,X)$
and the fact that $S(x_{(-1)})x_{(0)}=x$, we get
\begin{align*}
\beta_0\alpha_0 (h^0\otimes x)
 = & \left(S(h^0_{(3)})h^0_{(4)}h^0_{(1)}\otimes S(h^0_{(2)})x\right)
 = (h^0_{(1)}\otimes S(h^0_{(2)})x)
 = (h^0\otimes S(x_{(-1)})x_{(0)})
 = (h^0\otimes x)
\end{align*}
For $n\geq 1$ and for $(h^0\otimes\cdots\otimes h^n\otimes x)$ from
$\B{CM}_n(H,X)$ one has
\begin{align*}
\alpha_n
   & (h^0\otimes\cdots\otimes h^n\otimes x) \\
 = & \left(h^0_{(1)}\otimes h^0_{(1)}h^1_{(2)}\otimes\cdots\otimes
           h^0_{(n+1)}h^1_{(n)}\cdots h^{n-1}_{(2)}h^n_{(1)}\otimes 
           h^0_{(n+2)}h^1_{(n+1)}\cdots h^{n-1}_{(3)}h^n_{(2)}x\right)\\
 = & \left(h^0_{(1)}\epsilon(h^1_{(1)}\cdots h^n_{(1)})
           \otimes h^0_{(2)}h^1_{(2)}\epsilon(h^2_{(2)}\cdots h^n_{(2)})
           \otimes\cdots\otimes 
           h^0_{(n)}\cdots h^{n-1}_{(n)}\epsilon(h^n_{(n)})\right.\\
   & \hspace{1.5cm}\left.\otimes
           h^0_{(n+1)}\cdots h^n_{(n+1)}\otimes 
           h^0_{(n+2)}\cdots h^n_{(n+2)}x\right)\\
 = & \left(h^0_{(1)(1)}\epsilon(h^1_{(1)(1)}\cdots h^n_{(1)(1)})
           \otimes h^0_{(1)(2)}h^1_{(1)(2)}\epsilon(h^2_{(1)(2)}
           \cdots h^n_{(1)(2)})\otimes\cdots\otimes 
           h^0_{(1)(n)}\cdots h^{n-1}_{(1)(n)}\epsilon(h^n_{(1)(n)})
           \otimes\right.\\
   & \hspace{1.5cm}\left.\otimes
           h^0_{(2)}\cdots h^n_{(2)}\otimes 
           h^0_{(3)}\cdots h^n_{(3)}x\right)
\end{align*}
However, since $(h^0\otimes\cdots\otimes h^n\otimes x)$ is in
$\B{CM}^a_*(H,X)$, we have 
\begin{align*}
(h^0_{(1)}\otimes\cdots\otimes h^n_{(1)}\otimes 
 h^0_{(2)}\cdots h^n_{(2)}\otimes h^0_{(3)}\cdots h^n_{(3)}\otimes x)
 = & (h^0\otimes\cdots\otimes h^n\otimes 
      x_{(-2)}\otimes x_{(-1)}\otimes x_{(0)})
\end{align*}
This implies
\begin{align*}
\alpha_n
   & (h^0\otimes\cdots\otimes h^n\otimes x) \\
 = & \left(h^0_{(1)}\epsilon(h^1_{(1)}\cdots h^n_{(1)})
           \otimes h^0_{(2)}h^1_{(2)}\epsilon(h^2_{(2)}
           \cdots h^n_{(2)})\otimes\cdots\otimes 
           h^0_{(n)}\cdots h^{n-1}_{(n)}\epsilon(h^n_{(n)})
           \otimes x_{(-2)} \otimes x_{(-1)}x_{(0)}\right)\\
 = & \left(h^0_{(1)}\epsilon(h^1_{(1)}\cdots h^n_{(1)})
           \otimes h^0_{(2)}h^1_{(2)}\epsilon(h^2_{(2)}
           \cdots h^n_{(2)})\otimes\cdots\otimes 
           h^0_{(n)}\cdots h^{n-1}_{(n)}\epsilon(h^n_{(n)})
           \otimes x_{(-1)} \otimes x_{(0)}\right)
\end{align*}
by using the fact that $x_{(-1)}x_{(0)}=x$ for all $x\in X$.  Then
\begin{align*}
\beta_n\alpha_n
   & (h^0\otimes\cdots\otimes h^n\otimes x) \\
 = & \beta_n\left(h^0_{(1)}\epsilon(h^1_{(1)}\cdots h^n_{(1)})
           \otimes h^0_{(2)}h^1_{(2)}\epsilon(h^2_{(2)}
           \cdots h^n_{(2)})\otimes\cdots\otimes 
           h^0_{(n)}\cdots h^{n-1}_{(n)}\epsilon(h^n_{(n)})
           \otimes x_{(-1)} \otimes x_{(0)}\right)\\
 = & \left(S(x_{(-1)(3)})x_{(0)(-1)}h^0_{(1)}
           \epsilon(h^1_{(1)}\cdots h^n_{(1)})\otimes
           h^1_{(2)}\epsilon(h^2_{(2)}\cdots h^n_{(2)})\otimes
           \cdots\otimes x_{(-1)(1)}\otimes 
           S(x_{(-1)(2)})x_{(0)(0)}\right)\\
 = & \left(h^0_{(1)}\epsilon(h^1_{(1)}\cdots h^n_{(1)})\otimes
           h^1_{(2)}\epsilon(h^2_{(2)}\cdots h^n_{(2)})\otimes
           \cdots\otimes h^{n-1}_{(n)}\epsilon(h^n)
           \otimes x_{(-1)}\otimes x_{(0)}\right)\\
 = & \left(h^0_{(1)}\epsilon(h^1_{(1)}\cdots h^n_{(1)})\otimes
           h^1_{(2)}\epsilon(h^2_{(2)}\cdots h^n_{(2)})\otimes
           \cdots\otimes h^{n-1}_{(n)}\epsilon(h^n_{(n)})
           \otimes h^1_{(n+1)}\cdots h^n_{(n+1)}\otimes x\right)\\
 = & (h^0\otimes\cdots\otimes h^n\otimes x)
\end{align*}
as we wanted to show.
\end{proof}

\begin{lem}\label{graded_monomorphism}
Let $\B{T}^c_*(H,X)^\vee\xra{q_*}\B{CM}^c_*(H,X)^\vee :={}_H
\B{T}^c_*(H,X)^\vee$ be the quotient map.  Then
$q_*\alpha_*\beta_*=q_*$.
\end{lem}

\begin{proof}
The proof will be by direct calculation. For $n=0$ and $(h^0\otimes x)$
from $\B{T}^c_0(H,X)$ consider
\begin{align*}
\alpha_0\beta_0(h^0\otimes x)
 = & \alpha_0\left(S(h^0_{(3)})x_{(-1)}h^0_{(1)}\otimes 
                   S(h^0_{(2)})x_{(0)}\right)\\
 = & \left(S(h^0_{(3)(2)})x_{(-1)(1)}h^0_{(1)(1)}\otimes
           S(h^0_{(3)(1)})x_{(-1)(2)}h^0_{(1)(2)}S(h^0_{(2)})x_{(0)}
     \right)\\
 = & \left(S(h^0_{(3)})x_{(-2)}h^0_{(1)}\otimes 
           S(h^0_{(2)})x_{(-1)}x_{(0)}\right)\\
 = & S(h^0_{(2)})x_{(-1)}\cdot\left(h^0_{(1)}\otimes x_{(0)}\right)
\end{align*}
Then, 
\begin{align*}
q_0\alpha_0\beta_0(h^0\otimes x)
 = & \epsilon(S(h^0_{(2)})x_{(-1)})
     q_0\left(h^0_{(1)}\otimes x_{(0)}\right)\\
 = & q_0(h^0\otimes x)
\end{align*}
as we wanted to show.  Let $n\geq 0$ and $(h^0\otimes\cdots\otimes
h^{n+1}\otimes x)$ be from $\B{T}^c_*(H,X)^\vee$.  Consider
\begin{align*}
\alpha_{n+1}\beta_{n+1}
   & (h^0\otimes\cdots\otimes h^{n+1}\otimes x)\\
 = & \alpha_n\left(S(h^{n+1}_{(3)})x_{(-1)}h^0_{(1)}\otimes
           S(h^0_{(2)})h^1_{(1)}\otimes\cdots\otimes
           S(h^n_{(2)})h^{n+1}_{(1)}\otimes 
           S(h^{n+1}_{(2)})x_{(0)}\right)\\
 = & S(h^{n+1}_{(3)})x_{(-1)}\cdot
           \left(h^0_{(1)(1)}\otimes\alpha_n\left(
           h^0_{(1)(2)}S(h^0_{(2)})h^1_{(1)}\otimes\cdots\otimes
           S(h^n_{(2)})h^{n+1}_{(1)}\otimes 
           S(h^{n+1}_{(2)})x_{(0)}\right)\right)\\
 = & S(h^{n+1}_{(3)})x_{(-1)}\cdot
           \left(h^0\otimes\alpha_n\left(
           h^1_{(1)}\otimes\cdots\otimes S(h^n_{(2)})h^{n+1}_{(1)}
           \otimes S(h^{n+1}_{(2)})x_{(0)}\right)\right)
\end{align*}
However,
\begin{align*}
\alpha_n\left(z^0_{(1)}\otimes S(z^0_{(2)})z^1_{(1)}\otimes\cdots\otimes
         S(z^{n-1}_{(2)})z^n_{(1)}\otimes S(z^n_{(2)})x\right)
 = & (z^0\otimes\cdots\otimes z^n\otimes x)
\end{align*}
for any $(z^0\otimes\cdots\otimes z^n\otimes x)$ from $\B{T}^c_n(H,X)$.
Therefore
\begin{align*}
\alpha_{n+1}\beta_{n+1}(h^0\otimes\cdots\otimes h^{n+1}\otimes x)
 = & S(h^{n+1}_{(2)})x_{(-1)}\cdot
      (h^0\otimes\cdots\otimes h^{n+1}_{(1)}\otimes x_{(0)})
\end{align*}
Now apply $q_{n+1}$ on both sides to get the result.
\end{proof}

\begin{lem}\label{factorization}
$\B{T}^c_*(H,X)^\vee\xra{\beta_*}\B{T}^a_*(H,X)$ factors as
\begin{align*}
\B{T}^c_*(H,X)^\vee\xra{q_*}\B{CM}^c_*(H,X)^\vee
\xra{\beta'_*}\B{T}^a_*(H,X)
\end{align*}
iff $X$ is a stable anti-Yetter-Drinfeld module.
\end{lem}

\begin{proof}
Assume $X$ is a stable anti-Yetter-Drinfeld module.  Since
$\B{CM}^c_*(H,X)^\vee$ is defined as ${}_H\B{T}^c_*(H,X)^\vee$, we must
show that $\beta_*L_h=\epsilon(h)\beta_*$ for all $h\in H$.  Therefore,
consider
\begin{align*}
\beta_n
   & \left(h\cdot(h^0\otimes\cdots\otimes h^n\otimes x)\right)\\
 = & \beta_n(h_{(1)}h^0\otimes\cdots\otimes h_{(n+1)}h^n\otimes 
             h_{(n+2)}x)\\
 = & \left(S(h^n_{(3)})S(h_{(n+1)(3)})
           h_{(n+2)(1)}x_{(-1)}S^{-1}(h_{(n+2)(3)}) h_{(1)(1)}h^0_{(1)}
           \otimes S(h^1_{(2)})S(h_{(1)(2)})h_{(2)(1)}h^1_{(1)}
           \otimes\cdots\right.\\
   & \hspace{1.5cm}\left.
            \otimes S(h^{n-1}_{(2)})S(h_{(n)(2)})h_{(n+1)(1)}h^n_{(1)}
            \otimes S(h^n_{(2)})S(h_{(n+1)(2)})h_{(n+2)(2)}x_{(0)}
            \right)\\
 = & \epsilon(h)\left(S(h^n_{(3)})x_{(-1)}h^0_{(1)}\otimes
           S(h^0_{(2)})h^1_{(1)}\otimes\cdots\otimes
           S(h^{n-1}_{(2)})h^n_{(1)}\otimes S(h^n_{(2)})x_{(0)}\right)
\end{align*}
as we wanted to show.  On the opposite direction, assume $\beta$ factors
as $\beta_*=\beta'_*q_*$, which is to say $\beta_* L_h = \epsilon(h)
\beta_*$ for all $h\in H$.  Then
\begin{align*}
\left(S(h)x\right)_{(-1)}\otimes \left(S(h)x\right)_{(0)}
 = & \beta_0(\B{I}\otimes S(h)x)\\
 = & \beta_0\left(S(h_{(2)})h_{(3)}\otimes S(h_{(1)})x\right)\\
 = & \beta_0\left(S(h_{(1)})\cdot(h_{(2)}\otimes x)\right)\\
 = & \epsilon(S(h_{(1)})\beta_0(h_{(2)}\otimes x)\\
 = & \beta_0(h\otimes x)\\
 = & S(h_{(3)})x_{(-1)}h_{(1)}\otimes S(h_{(2)})x_{(0)}
\end{align*}
This finishes the proof.
\end{proof}

\begin{thm}
Let $X$ be a stable anti-Yetter-Drinfeld module.  Then
$\B{CM}^a_*(H,X)\xra{q_*\alpha_*}\B{CM}^c_*(H,X)^\vee$ is an isomorphism
of cyclic modules.
\end{thm}

\begin{proof}
By Lemma~\ref{factorization}, $\beta_*$ factors as $\beta_*=\beta'_*q_*$
as a morphism of para-cyclic modules.  Therefore the result we obtained
in Lemma~\ref{graded_epimorphism} reads as $\beta_*\alpha_* =
\beta'_*q_*\alpha_* = id_*$ restricted to $\B{CM}^a_*(H,X)$.  On the
other hand the result in Lemma~\ref{graded_monomorphism} reads as $q_* =
q_*\alpha_*\beta'_*q_*$.  Since $q_*$ is an epimorphism, every element
$\xi$ of $\B{CM}^c_*(H,X)$ is of the form $\xi = q_*({\bf h}\otimes x)$
for some $({\bf h}\otimes x)$ from $\B{T}^c_*(H,X)^\vee$.  Therefore
\begin{align*}
\xi = q_*({\bf h}\otimes x) 
    = q_*\alpha_*\beta'_*q_*({\bf h}\otimes x)
    = q_*\alpha_*\beta'_*(\xi)
\end{align*}
which proves $q_*\alpha_*$ and $\beta'_*$ are inverses of each other.
Since $\beta'_*$ is a morphism of para-cyclic modules, $q_*\alpha_*$
becomes an isomorphism of cyclic modules.
\end{proof}

\section{Computations}\label{Computations}

In order to simplify the computations, we assume $k=\B{C}$ in this
section.

\begin{exm}
Let $G$ be a discrete group and let $H=k[G]$ be the Hopf algebra of the
group ring of $G$ over $k$.  Consider $k$ as a trivial $G$--module via
$k[G]\xra{\epsilon}k$ and as a $k[G]$--comodule via the trivial coaction
$\rho_k(1) = 1\otimes 1$.  Then one can easily see that $k$ is a stable
$H$--module/comodule.  Moreover, $\B{T}^a_*(k[G],k)$ is $CC_*(k[G])$ the
classical cyclic object associated to the associative algebra $k[G]$
\cite[6.1.12]{Loday:CyclicHomology}.  Then $\B{CM}^a_*(k[G],k)$ consists
of elements of the form $\sum_i c_i(g^0_i\otimes\cdots\otimes
g^n_i\otimes 1)$ which satisfy
\begin{align*}
\sum_i c_i(g^0_i\otimes\cdots\otimes g^n_i\otimes 1\otimes 1)
 = & \sum_i c_i(g^0_i\otimes\cdots\otimes g^n_i\otimes 1\otimes
                g^0_i\cdots g^n_i)
\end{align*}
Since $k[G]^{\otimes n}$ is free over $k$ with basis from $G^{\times
n+1}$, this implies $g^0_i\cdots g^n_i = 1$ for any
$(g^0_i\otimes\cdots\otimes g^n_i\otimes 1)$ in the summation $\sum_i
c_i(g^0_i\otimes\cdots\otimes g^n_i\otimes 1)$.  In other words,
$\B{CM}^a_*(k[G],k)$ is the $\left<1\right>$-component of $CC_*(k[G])$
which is denoted by $CC_*(k[G])_{\left<1\right>}$.  Then 
\begin{align}
HC^{\B{CM},a}_n (k[G],k) 
 := HC_n\B{CM}^a_*(k[G],k) 
 = HC_nCC_*(k[G])_{\left<1\right>} 
\end{align}
And, according to \cite{Loday:CyclicHomology}, one has
\begin{align}
HC_nCC_*(k[G])_{\left<1\right>} 
 \cong \bigoplus_{i\geq 0}H_{n-2i}(G)
\end{align}
for any $\geq 0$.
\end{exm}

\begin{exm}
Let $\G{g}$ be any Lie algebra and let $H=U(\G{g})$ be its universal
enveloping algebra.  Again, consider $k$ as a trivial $H$--comodule via
$1$.  Fix a character $U(\G{g})\xra{\delta}k$ and consider $k$ as a
$U(\G{g})$--module via this character.  Denote this one dimensional
stable $U(\G{g})$--module/comodule by $k_{(1,\delta)}$.  By using a
fixed basis for $\G{g}$, Lemma~\ref{main_lemma} and
Poincar\'e--Birkhoff--Witt Theorem one can conclude that
\begin{align*}
U(\G{g})^{\otimes n+1}\fotimes{U(\G{g})}k_{(1,\delta)} 
\cong \left(U(\G{g})^{\otimes n+1}\right)^{U(\G{g})} \cong k
\end{align*}
which implies $\B{CM}^a_*(U(\G{g}),k_{(1,\delta)})=CC_*(k)$.  This means
\begin{align}
HC^{\B{CM},a}_n(U(\G{g}),k_{(1,\delta)})
 := HC_n\B{CM}^a_*(U(\G{g}),k_{(1,\delta)}) \cong HC_n(k)
\end{align}
for any $n\geq 0$.
\end{exm}

\begin{exm}\label{quantum}
Let $\G{g}$ be a semi-simple Lie algebra of rank $N$ and let
$U_q(\G{g})$ be the quantum deformation of the Lie algebra $\G{g}$.  One
can recall the presentation of $U_q(\G{g})$ from
\cite{Kaygun:BialgebraCyclic}.  Fix a group-like element
$K_I=K_1^{a_1}\cdots K_N^{a_N}$ from $U_q(\G{g})$ where $a_i\in\B{Z}$
and $I=\left<a_1,\ldots,a_N\right>$.  Consider $k_I=k$ as a
$U_q(\G{g})$--module/comodule via the counit $\epsilon$ and the
grouplike element $K_I$.  One can see that $k_I$ is a stable
$U_q(\G{g})$--module/comodule.  Because of Lemma~\ref{main_lemma} and
the quantum Poincar\'e--Birkhoff--Witt Theorem
\begin{align*}
U_q(\G{g})^{\otimes n+1}\fotimes{U_q(\G{g})}k_I
\cong \left(U_q(\G{g})^{\otimes n+1}\otimes k_I\right)^{U_q(\G{g})}
\cong k\left[K_1^\pm,\ldots,K_N^\pm\right]_{\left<K_I\right>}
\end{align*}
This means $\B{CM}^a_*(U_q(\G{g}),k_I)\cong
CC_*(k\left[K_1^\pm,\ldots,K_N^\pm\right])_{\left<K_I\right>}$ which in
turn is isomorphic to the cyclic object
$CC_*(k\left[K_1^\pm,\ldots,K_N^\pm\right])_{\left<1\right>}$ since the
group $\left<K_1^\pm,\ldots,K_N^\pm\right>$ is abelian.  Then
\begin{align}
HC^{\B{CM},a}_n(U_q(\G{g}),k_I) 
 \cong & HC_nCC_*(k[\B{Z}^{\times N}])_{\left<1\right>} \cong
 \bigoplus_{i\geq 0}H_{n-2i}(\B{Z}^{\times N})
\end{align}
which implies
\begin{align}
HC^{\B{CM},a}_n(U_q(\G{g}),k_I) 
 \cong & \begin{cases}
         k            & \text{  if $n$ is even }\\
         k^{\oplus N} & \text{  if $n$ is odd }
	 \end{cases}
\end{align}
\end{exm}

\begin{exm}
Let $\C{H}(N)$ be the Hopf algebra of codimension $N$ foliations.  One
can recall the presentation from \cite{Kaygun:BialgebraCyclic}.  Again,
fix a character $\C{H}(N)\xra{\delta}k$ and use $k$ as a stable
$\C{H}(N)$--module/comodule via the pair $(1,\delta)$.  Since there are
no group-like elements in $\C{H}(N)$ except 1,
\begin{align}
\C{H}(N)^{\otimes n+1}\fotimes{\C{H}(N)}k_{(1,\delta)}
\cong \left(\C{H}(N)^{\otimes n+1}\otimes k_{(1,\delta)}\right)^{\C{H}(N)}
\cong k
\end{align}
which implies $\B{CM}^a_*(\C{H}(n),k_{(1,\delta)})=CC_*(k)$.  Therefore
\begin{align}
HC^{\B{CM},a}_n(\C{H}(N),k_{(1,\delta)})
 := HC_n\B{CM}^a_*(\C{H}(N),k_{(1,\delta)}) \cong HC_n(k)
\end{align}
for any $n\geq 0$.  This result is in direct contract with the dual
theory.  In \cite{Tamas:Thesis} T. Antal proved that if one takes the
character $\delta$ which satisfies
\begin{align}
\delta(X) = &\ 0 = \delta(\delta_n)  & \delta(Y) = 1
\end{align}
for any $n\geq 1$, then $(1,\delta)$ is a modular pair and the classical
Hopf cyclic homology group $HC^{(1,\delta)}_1(\C{H}(1))$, which is the
same as $HC^{\B{CM},c}_1(\C{H}(1),k_{(1,\delta)})$, is two dimensional.
Hence $HC^{\B{CM},c}_*(\C{H}(1),k_{(1,\delta)})$ is different than
$HC_*(k)$.
\end{exm}

\bibliographystyle{plain}
\bibliography{bibliography}

\vspace{2cm}
\noindent{\small\sc 
Department of Mathematics, Ohio State University,
Columbus, Ohio 43210, USA}

\noindent{\it E-mail address:}\ {\tt kaygun@math.ohio-state.edu}

\end{document}